\documentclass[11pt,a4paper]{article}

\usepackage[T1]{fontenc}
\usepackage{epsfig}
\usepackage{amsmath, amsthm}
\usepackage{amsfonts,amssymb}

\newtheorem{conj}{Conjecture}[section]
\newtheorem{theo}[conj]{Theorem}

\newtheorem{lem}[conj]{Lemma}
\newtheorem{prop}[conj]{Proposition}
\newtheorem{coro}[conj]{Corollary}

\newcommand{\de}{\mathrm{d}}
\newcommand{\e}{\mathrm{e}}
\newcommand{\conv}{\mathrm{conv}}

\newcommand{\interior}{\mathrm{int}}

\newcommand{\supp}{\mathrm{supp}}
\newcommand{\dist}{\mathrm{dist}}
\newcommand{\Var}{\mathrm{Var}}
\newcommand{\Hess}{\mathrm{Hess}}

\newcommand{\R}{\mathbb{R}}

\newcommand{\eps}{\varepsilon}

\begin{document}

\title{Concavity properties of extensions of the parallel volume}

\date{}

\author{Arnaud Marsiglietti}
\maketitle
{\raggedright Universit\'e Paris-Est, LAMA (UMR 8050), UPEMLV, UPEC, CNRS, F-77454, Marne-la-Vall\'ee, France} \\

\noindent
arnaud.marsiglietti@u-pem.fr \\

{\raggedright {\scriptsize The author was supported in part by the Agence Nationale de la Recherche, project GeMeCoD (ANR 2011 BS01 007 01).}}

\begin{abstract}

In this paper we establish concavity properties of two extensions of the classical notion of the outer parallel volume. On the one hand, we replace the Lebesgue measure by more general measures. On the other hand, we consider a functional version of the outer parallel sets.

\end{abstract}


\section{Introduction}

As an analogue of the famous {\it concavity of entropy power} in Information theory (see $e.g.$ \cite{C}, \cite{V}), Costa and Cover \cite{Costa} conjectured that the $\frac{1}{n}$-power of the {\it parallel volume} $|A+tB_2^n|$ is a concave function, where $|\cdot|$ denotes the Lebesgue measure and where $B_2^n$ denotes the Euclidean closed unit ball.

\begin{conj}[Costa-Cover \cite{Costa}]\label{theo:1.1}

Let $A$ be a bounded measurable set in $\mathbb{R}^n$ then the function $ t \mapsto |A + t B_2^n|^{\frac{1}{n}} $ is concave on $\mathbb{R}_+$.

\end{conj}

This conjecture has been studied by the author and Fradelizi in \cite{nous}, where it was shown to be true for any sets in dimension $1$ and for any connected sets in dimension $2$. However they showed that the conjecture fails for arbitrary sets in dimension $2$ and for arbitrary connected sets in dimension greater than or equal to $3$. \\

The notion of parallel volume can be extended by considering measures that are more general than the Lebesgue measure. An extension, provided by Borell in \cite{Borell1}, \cite{Borell}, follows from the Brunn-Minkowski inequality, which states that for every $\lambda \in [0,1]$ and for every compact subsets $A,B$ of $\R^n$,
\begin{eqnarray}\label{BM}
|(1-\lambda)A + \lambda B|^{\frac{1}{n}} \geq (1-\lambda)|A|^{\frac{1}{n}} + \lambda |B|^{\frac{1}{n}}.
\end{eqnarray}
Borell defined for $s \in [- \infty, + \infty]$ the \textit{$s$-concave measures}, satisfying by definition a similar inequality to (\ref{BM}):
$$ \mu((1 - \lambda)A + \lambda B) \geq M_s^{\lambda}(\mu(A),\mu(B)) $$
for every $\lambda \in [0,1]$ and for every compact subsets $A,B$ of $\mathbb{R}^n$ such that $\mu(A)\mu(B)>0$, where $M_s^{\lambda}(a,b)$ denotes the $s$-mean of the non-negative real numbers $a, b$ with weights $1-\lambda$ and $\lambda$, defined as
$$ M_s^{\lambda}(a,b) = ((1-\lambda)a^s + \lambda b^s)^{\frac{1}{s}} \quad \mbox{if $s \notin \{-\infty, 0, +\infty\}$}, $$
$M_{-\infty}^{\lambda}(a,b) = \min(a,b)$, $M_0^{\lambda}(a,b) = a^{1-\lambda} b^{\lambda}$, $M_{+\infty}^{\lambda}(a,b) = \max(a,b)$.

Basic properties of the $s$-concave measures are given in the next section. The case $s=0$ corresponds to log-concave measures. The most famous example of a log-concave measure is the {\it standard multivariate Gaussian measure}
$$ \de \gamma_n(x) = \frac{1}{(2 \pi)^{\frac{n}{2}}} \e^{-\frac{|x|^2}{2}} \de x, $$
where $|\cdot|$ stands for the Euclidean norm. These measures are of interest. For example, isoperimetric inequalities have been established for the Gaussian measure $\gamma_n$ by Borell \cite{Borell0} and independently by Sudakov and Cirel'son \cite{Sudakov}, which state that among sets of given Gauss measure, half-spaces minimize the Gauss surface area (see also \cite{KLS}, \cite{Bobkov4}, \cite{Cianchi}, \cite{Castro}).

In this paper, we pursue the study of these measures by considering the following problem: \\

\noindent
{\bf Problem 1.}
{\it
Let $s \in [-\infty, +\infty]$. Let $\mu$ be an $s$-concave measure on $\R^n$ and let $A$ be a compact subset of $\mathbb{R}^n$. Is the function $t \mapsto \mu(A + t B_2^n)$ $s$-concave on $\mathbb{R}_+$?
} \\

Notice that from the Brunn-Minkowski inequality~(\ref{BM}), the $n$-dimensional Lebesgue measure is $\frac{1}{n}$-concave. Thus Problem~1 generalizes the Costa-Cover conjecture. \\

Other extensions of geometric inequalities can be set up by considering functional versions. The most famous extension of this type in the Brunn-Minkowski theory is certainly the Pr\'ekopa-Leindler inequality (see \cite{Leindler}, \cite{Prekopa}, \cite{Borell}). Functional versions provide new proofs of geometric inequalities as well as new applications. Other examples of such extensions are a functional version of the Blaschke-Santal\'o inequality and a functional version of the Mahler conjecture (see $e.g.$ \cite{Ball}, \cite{AKM}, \cite{FradeliziM1}, \cite{FradeliziM2}, \cite{Lehec}).

To establish a functional version of Conjecture~\ref{theo:1.1}, we consider a functional version of {\it parallel sets} $A + t B_2^n$. We set up the following problem (the notion of $\gamma$-concave functions is defined in the next section): \\

\noindent
{\bf Problem 2.}
{\it
Let $\gamma \geq -\frac{1}{n}$. Let $f: \R^n \to \R_+$ be a measurable non-negative function and $g : \R^n \to \R_+$ be a $\gamma$-concave function. Let us denote
$$ h_t^{(\gamma)}(z) = \underset{f(x)>0;~g(y)>0}{\sup_{z=x+ty}}(f(x)^{\gamma} + tg(y)^{\gamma})^{\frac{1}{\gamma}} \quad \mbox{and} \quad h_t^{(0)}(z) = \sup_{z=x+ty}f(x)g(y)^t. $$
Is the function $ t \mapsto \int_{\R^n} h_t^{(\gamma)}(z) \, \de z $ $\frac{\gamma}{1+\gamma n}$-concave on $\R_+$?
} \\

The Costa-Cover conjecture is a particular case of Problem~2 by taking $f=1_A$, $g=1_{B_2^n}$ and $\gamma = 0$. 

For $\gamma < 0$, one can connect the function $h_t^{(\gamma)}$ with the Hopf-Lax solution of the Hamilton-Jacobi equation (see Section 4 for precise definitions). The Hopf-Lax solution is of interest. For example, it can be used to show that {\it hypercontractivity} of this solution is equivalent to {\it log-Sobolev inequalities} (see $e.g.$ \cite{Bobkov3}, \cite{Gozlan}). In Problem~2, we pursue the study of this solution by investigating concavity properties in time of the Hopf-Lax solution of the Hamilton-Jacobi equation. \\

We shall prove that both Problem~1 and Problem~2 have positive answers in dimension $1$. However, the Costa-Cover conjecture fails in dimension $n \geq 2$ in such a generality, thus we won't expect positive answers of these more general problems in dimension $n \geq 2$. Using the geometric localization theorem of Kannan, Lov\'asz and Simonovits \cite{KLS} in the form established by Fradelizi and Gu\'edon \cite{Fradelizi} (see Theorem~\ref{FG} below), we prove: \\

\noindent
{\bf Theorem 1.}
{\it
Let $s \in [- \infty, \frac{1}{2}] \cup [1, + \infty]$. Let $\mu$ be an $s$-concave measure on $\R$ and let $A$ be a compact subset of $\R$. Then the function $t \mapsto \mu(A + t [-1,1])$ is $s$-concave on $\mathbb{R}_+$.
} \\

Moreover, for $s \in (\frac{1}{2}, 1)$ there exists a compact subset $A$ of $\R$ such that $t \mapsto \mu(A + t [-1,1])$ is not $s$-concave on $\mathbb{R}_+$ (see Proposition~\ref{gap} below).

Using a precise analysis of the Hopf-Lax solution, we prove: \\

\noindent
{\bf Theorem 2.}
{\it
Let $\gamma \in [-1, 0]$. Let $f: \R \to \R_+$ be such that $f^{\gamma}$ (to be interpreted by $- \log(f)$ when $\gamma = 0$) is a bounded Lipschitz function. Define for every $y \in \R$, $V(y) =  \frac{|y|^p}{p}$ with $p \geq 1$, and
$$ h_t^{(\gamma)}(z) = \underset{f(x)>0;~V(y)>0}{\sup_{z=x+ty}}(f(x)^{\gamma} + tV(y))^{\frac{1}{\gamma}}, \qquad h_t^{(0)}(z) = \sup_{z=x+ty}f(x)\e^{-t V(y)}. $$
If there exists $t_0 > 0$ such that for almost every $z \in \R$, $t \mapsto h_t^{(\gamma)}(z)$ is twice differentiable on $[0,t_0]$ with derivatives that are bounded by an integrable function and if $\lim_{z \to \pm \infty}\frac{\partial h_t^{(\gamma)}}{\partial z} (z) = 0$, then the function $ t \mapsto \int_{\R} h_t^{(\gamma)}(z) \, \de z $
is concave on $[0,t_0]$.
} \\

Notice that the conclusion of Theorem~2 is stronger than the expected conclusion in Problem~2 (when restricted to log-concave functions of the form $|y|^p/p$, $p \geq 1$). \\

In the next section, we present some notation and we explain basic properties of the $s$-concave measures. In section 3, we give a complete answer to Problem~1. In section 4, we give a partial answer to Problem~2. To conclude this paper, we derive a weighted Brascamp-Lieb-type inequality from our functional version.

\section{Preliminaries}

We work in the Euclidean space $\mathbb{R}^n$, $n \geq 1$, equipped with the usual scalar product $\langle \cdot,\cdot \rangle$ and the $\ell_2^n$ norm $|\cdot |$. The closed unit ball is denoted by $B_2^n$ and the canonical basis by $(e_1, \cdots, e_n)$. We also denote by $|\cdot |$ the Lebesgue measure on $\mathbb{R}^n$. For non-empty sets $A,B \subset \mathbb{R}^n$ we define the \textit{Minkowski sum}
$$ A+B = \{a + b ; \, a \in A, b \in B \}.$$
We denote by $\interior(A)$ the \textit{interior} of the set $A$, by $\overline{A}$ the {\it closure} of $A$ and by $\partial A = \overline{A} \setminus \interior(A)$ the {\it boundary} of $A$.

For an arbitrary (non-negative) measure $\mu$, we call {\it outer parallel $\mu$-volume} of a compact set $A$ the function defined on $\R_+$ by $t \mapsto V_A^{\mu}(t)=\mu(A + t B_2^n)$. We simply call this function {\it outer parallel volume} when $\mu$ is the Lebesgue measure.

Let us recall basic properties of the $s$-concave measures. For $s \leq \frac{1}{n}$, Borell showed in \cite{Borell} that any $s$-concave measure which is absolutely continuous with respect to the Lebesgue measure on $\R^n$ admits a $\gamma$-concave density function, with
$$ \gamma = \frac{s}{1-sn}  \in [-\frac{1}{n}, + \infty], $$
where a function $f$ is said to be \textit{$\gamma$-concave}, with $\gamma \in [-\infty, +\infty]$, if the inequality
$$ f((1 - \lambda)x + \lambda y) \geq M_{\gamma}^{\lambda}(f(x), f(y)) $$
holds for every $\lambda \in [0,1]$ and for every $x,y \in \R^n$ such that $f(x)f(y)>0$. Notice that a 1-concave function is concave on its support, that the level sets of a $- \infty$-concave function are convex, and that a $+ \infty$-concave function is constant on its support. For $s > 1$, Borell showed that every (non-null) $s$-concave measure is a Dirac measure. Notice that an $s$-concave measure is $r$-concave for all $r \leq s$.

\section{The $s$-concavity of the parallel $\mu$-volume}

In this section we solve Problem~1.

Let us note that since the function $t \mapsto \mu(A + t B_2^n)$ is non-decreasing, it follows that the answer to Problem~1 is positive for $s=-\infty$ (for every non-negative measures). Note also that Problem~1 is solved for convex sets. Indeed, if $A$ is a compact convex subset of $\R^n$, then for every $\lambda \in [0,1]$ and for every $t_1,t_2 \in \R_+$, one has
\begin{eqnarray*}
\mu(A + ((1-\lambda)t_1 + \lambda t_2) B_2^n) & = & \mu((1-\lambda)(A + t_1B_2^n) + \lambda (A + t_2B_2^n)) \\ & \geq & ((1-\lambda) \mu(A + t_1B_2^n)^s + \lambda \mu(A + t_2B_2^n)^s)^{\frac{1}{s}}.
\end{eqnarray*}

Before proving Theorem~1, we establish a preliminary lemma in dimension $1$.

\begin{lem}\label{raccord}

Let $\mu$ be an $s$-concave measure on $\R$ which is absolutely continuous with respect to the Lebesgue measure and let $A$ be a compact subset of $\R$. Then the function $t \mapsto V_A^{\mu}(t) = \mu(A + t [-1,1])$ admits left and right derivatives on $(0, + \infty)$, denoted respectively by $(V_A^{\mu})'_-$ and by $(V_A^{\mu})'_+$, which satisfy for every $t > 0$, $(V_A^{\mu})'_-(t) \geq (V_A^{\mu})'_+(t)$.

\end{lem}

\begin{proof}
Let us denote by $\psi$ the density of the measure $\mu$ and let us denote by $[\alpha,\beta]$ the support of $\mu$, with $- \infty \leq \alpha < \beta \leq + \infty$. Notice that for every $t > 0$, the set $A + t [-1,1]$ is a finite disjoint union of intervals, thus one can assume that $A = \cup_{i=1}^N [a_i,b_i]$, with $\alpha \leq a_1 < b_1 < \dots < a_N < b_N \leq \beta$, and $a_1 = \alpha$ when $\alpha \in A$, $b_N = \beta$ when $\beta \in A$. Let us denote $t_i = \frac{a_{i+1} - b_i}{2}$, $i \in \{ 0, \cdots, N \}$, with the convention that $b_0 = 2 \alpha - a_1$ and that $a_{N+1} = 2 \beta - b_N$. Notice that $V_A^{\mu}$ is differentiable on $(0, + \infty) \setminus \{t_0, \cdots, t_{N} \}$ and that for every $t > 0$, one has
$$ (V_A^{\mu})'_+(t) = \sum_{a \in \partial(A + t [-1,1])} \psi(a) 1_{(\alpha,\beta)}(a), \quad (V_A^{\mu})'_-(t) = \sum_{a \in \partial(A + t (-1,1))} \psi(a) 1_{[\alpha,\beta]}(a). $$

Notice that $\overline{A + t[-1,1]} = \overline{A + t \, (-1,1)}$ thus $\partial(A + t [-1,1]) \subset \partial(A + t \,(-1,1))$. We conclude that for every $t > 0$, $(V_A^{\mu})'_-(t) \geq (V_A^{\mu})'_+(t)$.
\end{proof}

The proof of Theorem~1 uses the following localization theorem.

\begin{theo}[Fradelizi, Gu\'edon \cite{Fradelizi}]\label{FG}

Let $n$ be a positive integer, let $K$ be a compact convex set in $\R^n$ and denote by $\mathcal{P}(K)$ the set of probability measures on $\R^n$, with support contained in $K$. Let $f : K \to \R$ be an upper semi-continuous function, let $s \in [- \infty, \frac{1}{2}]$ and denote by $P_f$ the set of $s$-concave probability measures $\mu$, with support contained in $K$, satisfying $\int f \, \de \mu \geq 0$. Let $\Phi : \mathcal{P}(K) \to \R$ be a convex $w*$-upper semi-continuous function. Then $\sup\{ \Phi(\mu) ; \mu \in P_f \}$ is achieved at a probability measure $\nu$ which is either a Dirac measure at a point $x$ such that $f(x) \geq 0$, or a probability measure $\nu$ which is $s$-affine on a segment $[\alpha, \beta] \subset K$, such that $\int f \, \de \nu = 0$ and $\int_{\alpha}^x f \, \de \nu > 0$ on $(\alpha,\beta)$ or $\int_x^{\beta} f \, \de \nu > 0$ on $(\alpha,\beta)$.

\end{theo}

\begin{proof}[Proof of Theorem 1.]
Let $\mu$ be an $s$-concave measure on $\R$ and let $A$ be a compact subset of $\R$. 

If $s > 1$ then $\mu$ is a Dirac measure $\delta_x$, $x \in \mathbb{R}$, and one can see that the function $t \mapsto \delta_x(A + t [-1,1])$ is constant on its support. Thus this function is $+ \infty$-concave on $\mathbb{R}_+$, which proves Theorem~1 for $s>1$. Notice that this argument is also valid in higher dimension.

If $s=1$ then for every $x \in \R$, $\de \mu(x) = C 1_{[\alpha,\beta]}(x) \de x $, with $C>0$ is a constant and $[\alpha,\beta]$ is an interval of $\R$. One can see by a direct computation that $t \mapsto \mu(A + t[-1,1])$ is $1$-concave on $\R_+$.

If $s \leq \frac{1}{2}$ then $\mu$ admits a density with respect to the Lebesgue measure on $\R$ which is $\frac{s}{1-s}$-concave. We assume that $s \neq 0$, the case $s = 0$ follows by continuity. Let $t_0 > 0$ such that the set $A + t_0[-1,1]$ is an interval. Notice that the function $t \mapsto V_A^{\mu}(t) = \mu(A + t[-1,1])$ is $s$-concave on $[t_0, + \infty)$.

To prove that $V_A^{\mu}$ is $s$-concave on $(0,t_0)$, we shall reduce the problem to extremal measures $\nu$ described in Theorem~\ref{FG}. For these specific measures, the outer parallel volume is easy to compute and is twice differentiable except for a finite number of points. Then we prove by differentiation that the function $t \mapsto V_A^{\nu}(t) = \nu(A + t[-1,1])$ is $s$-concave outside points of non-differentiability. To conclude global $s$-concavity, we use Lemma~\ref{raccord}. \\

\noindent
{\bf Step 1:} \underline{Reduction to extremal measures} \\
Let $t_1,t_2 \in (0,t_0)$ such that $\mu(A + t_1 [-1,1])\mu(A + t_2 [-1,1])>0$. We assume that $t_1<t_2$ without loss of generality. We denote $K = A + t_0 [-1,1]$ and we consider the restriction of $\mu$ over $K$, thus it is a finite measure that we can assume to be a probability measure without loss of generality. For convenience, we still denote by $\mu$ this measure.  Let $0<\eps<t_2-t_1$. We apply Theorem~\ref{FG} to $f : K \to \R$ defined by 
$$ f = 1_{A + t_2 [-1,1]} - \tau_{\eps} 1_{A + (t_1+\eps) (-1,1) } $$
and $\Phi_{\eps} : \mathcal{P}(K) \to \R$ defined by
$$ \Phi_{\eps} = \rho_{\eps} 1_{A + t_1 [-1,1]} - 1_{A + \frac{(t_1+\eps) + t_2}{2} (-1,1) } $$
where
$$ \rho_{\eps} = \left( \frac{1}{2} \left( \frac{\mu(A + t_2 [-1,1])^s}{\mu(A + (t_1+\eps) [-1,1])^s} + 1 \right)  \right)^{\frac{1}{s}}, \quad \tau_{\eps} = \frac{\mu(A + t_2 [-1,1])}{\mu(A + (t_1+\eps) [-1,1])}. $$
The choice of introducing $\eps$ and the open interval $(-1,1)$ is a technical trick to get upper semi-continuous functions regarding $f$ and $\Phi_{\eps}$ and to make the argument work for Dirac measures in Step~2 below. We shall prove that $\Phi_{\eps}(\mu) \leq 0$, thus by letting $\eps$ go to $0$ and by using that $V_A^{\mu}$ is continuous on $(0,t_0)$, this will lead to the conclusion that $V_A^{\mu}$ is $s$-concave on $[t_1, t_2]$, for arbitrary $t_1,t_2 \in (0,t_0)$. To prove that $\Phi_{\eps}(\mu) \leq 0$, we shall prove that $\Phi_{\eps}(\nu) \leq 0$ for the extremal probability measures $\nu$ described in Theorem~\ref{FG}. First, notice that $\tau_{\eps} \geq 1$. If $\tau_{\eps} = 1$, then $V_A^{\mu}$ is constant on $[t_1+\eps, t_2]$. Thus $V_A^{\mu}$ is $s$-concave on $[t_1+\eps, t_2]$. Thereafter, we assume that $\tau_{\eps} > 1$. \\

\noindent
{\bf Step 2:} \underline{$s$-concavity for extremal measures} \\
- Let $\nu$ be a Dirac measure $\delta_x$ with $x$ such that $f(x) \geq 0$. The condition $f(x) \geq 0$ says that $1_{A + t_2 [-1,1]}(x) \geq \tau_{\eps} 1_{A + (t_1+\eps) (-1,1) }(x)$. Since $\tau_{\eps} > 1$, it follows that $x \notin A+(t_1+\eps) (-1,1)$. Thus $x \notin A+t_1[-1,1]$. Hence $\Phi_{\eps}(\delta_x) \leq 0$. \\

\noindent
- Let $\nu$ be an $s$-affine measure with support $[\alpha,\beta]$, $i.e.$ the density of the measure $\nu$, denoted by $\psi$, satisfies for every $x \in \R$, $\psi(x) = (mx + p)^{1/\gamma} 1_{[\alpha,\beta]}(x)$, where $\gamma = s/(1-s) \in [-1,1] \setminus \{0\}$ and where $m \in\R \setminus \{0\}$, $p \in \R$ are such that for every $x \in [\alpha,\beta]$, $mx+p \geq 0$. We assume that $\int_{[x,\beta]} f \, \de \nu < 0$ on $(\alpha,\beta)$, which means that for every $x \in (\alpha,\beta)$,
\begin{eqnarray}\label{condition}
\nu ((A+t_2 [-1,1]) \cap [x,\beta]) < \tau_{\eps} \nu ((A + (t_1+\eps) [-1,1]) \cap [x,\beta] ).
\end{eqnarray}
If $\beta \notin A+(t_1+\eps) [-1,1]$, then there exists $x \in (\alpha,\beta)$ such that $(A+(t_1+\eps) [-1,1]) \cap [x,\beta] = \emptyset$. This contradicts~(\ref{condition}). It follows that $\beta \in A+(t_1+\eps) [-1,1]$.

Notice that the function $V_A^{\nu}$ is twice differentiable on $[t_1 + \eps, t_2]$ outside a finite number of points $s_0, \dots, s_{\ell}$, with $s_0 := t_1 + \eps < \cdots < s_{\ell} < s_{\ell+1} := t_2$. \\

\noindent
{\bf Case 1:} $\alpha \in A+(t_1+\eps) [-1,1]$. \\
Let $j \in \{0, \cdots, \ell \}$. Notice that $A + s_j[-1,1]$ is a finite disjoint union of intervals containing $\alpha$ and $\beta$, hence we can assume that $A + s_j[-1,1] = \cup_{i=1}^{N} [a_i,b_i]$, with $a_1 \leq \alpha<b_1<\dots<a_N<\beta \leq b_N$. We assume that $N \geq 2$, otherwise the result clearly holds. Let $t \in [s_j, s_{j+1})$. We get
\begin{eqnarray*}
V_{A}^{\nu}(t) & = & \sum_{i=1}^N \int_{a_i - t}^{b_i + t} (mx + p)^{\frac{1}{\gamma}} 1_{[\alpha,\beta]}(x) \de x, \\ (V_{A}^{\nu})'(t) & = & \sum_{i=2}^N \left( (m(b_{i-1} + t) + p)^{\frac{1}{\gamma}} + (m(a_i - t) + p)^{\frac{1}{\gamma}} \right), \\ (V_{A}^{\nu})''(t) & = & \frac{m}{\gamma} \left( \sum_{i=2}^{N} \left( (m(b_{i-1} + t) + p)^{\frac{1-\gamma}{\gamma}} - (m(a_i - t) + p)^{\frac{1-\gamma}{\gamma}} \right) \right) \leq 0.
\end{eqnarray*}
Hence the function $V_{A}^{\nu}$ is concave on $[s_j, s_{j+1})$. We deduce that $V_{A}^{\nu}$ is piecewise $s$-concave on $[t_1 + \eps, t_2]$. From Lemma~\ref{raccord}, we get that for every $t \in [t_1 + \eps, t_2]$,
$$ s((V_A^{\nu})^s)'_-(t) = s^2 (V_A^{\nu})'_-(t) (V_A^{\nu})^{s-1}(t) \geq s^2 (V_A^{\nu})'_+(t) (V_A^{\nu})^{s-1}(t) = s((V_A^{\nu})^s)'_+(t). $$
We conclude that the function $V_{A}^{\nu}$ is $s$-concave on $[t_1 + \eps, t_2]$. \\

\noindent
{\bf Case 2:} $\alpha \notin A+(t_1+\eps) [-1,1]$. \\
Let $j \in \{0, \cdots, \ell \}$. If $\alpha \in A + s_j[-1,1]$ then from the previous case we can conclude that $V_{A}^{\nu}$ is $s$-concave on $[s_j, s_{j+1})$. Thus we can assume that $A + s_j[-1,1] = \cup_{i=1}^{N} [a_i,b_i]$, with $\alpha < a_1 <b_1<\dots<a_N<\beta\leq b_N $ and that $\alpha \notin A + s_{j+1}[-1,1]$. We assume that $N \geq 2$, otherwise the result clearly holds. In the following, we denote $a_i(t) = m(a_i - t) + p$, $b_i(t) = m(b_i + t) + p$, $1 \leq i \leq N-1$ and $a_N(t) = m(a_N - t) + p$, $b_N(t) = m\beta + p$.

Let $t \in [s_j, s_{j+1})$. We get
\begin{eqnarray*}
V_{A}^{\nu}(t) & = & \frac{1}{m} \frac{\gamma}{1 + \gamma} \sum_{i=1}^{N} \left( b_i(t)^{\frac{1+\gamma}{\gamma}} - a_i(t)^{\frac{1+\gamma}{\gamma}} \right), \\ (V_{A}^{\nu})'(t) & = & \sum_{i=1}^{N-1} \left( b_i(t)^{\frac{1}{\gamma}} + a_i(t)^{\frac{1}{\gamma}} \right) + a_N(t)^{\frac{1}{\gamma}}, \\ (V_{A}^{\nu})''(t) & = & \frac{m}{\gamma} \left( - a_1(t)^{\frac{1-\gamma}{\gamma}} + \sum_{i=2}^{N} \left( b_{i-1}(t)^{\frac{1-\gamma}{\gamma}} - a_i(t)^{\frac{1-\gamma}{\gamma}} \right) \right).
\end{eqnarray*}
Then the function $V_{A}^{\nu}$ is $s$-concave on $[s_j, s_{j+1})$ if and only if for every $t \in [s_j, s_{j+1})$, $V_{A}^{\nu}(t) (V_{A}^{\nu})''(t) \leq (1-s) (V_{A}^{\nu})'(t)^2$ if and only if
\begin{eqnarray*}
\left( \sum_{i=1}^{N} \left( a_i(t)^{\frac{1+\gamma}{\gamma}} - b_i(t)^{\frac{1+\gamma}{\gamma}} \right) \right) \left( \sum_{i=1}^{N-1} \left( a_i(t)^{\frac{1-\gamma}{\gamma}} - b_i(t)^{\frac{1-\gamma}{\gamma}} \right) + a_N(t)^{\frac{1-\gamma}{\gamma}} \right) \\ \leq \left( \sum_{i=1}^{N-1} \left( b_i(t)^{\frac{1}{\gamma}} + a_i(t)^{\frac{1}{\gamma}} \right) + a_N(t)^{\frac{1}{\gamma}} \right)^2.
\end{eqnarray*}

If $m \gamma > 0$, then one has 
$$ 0 \leq a_1(t)^{\frac{1-\gamma}{\gamma}} < b_1(t)^{\frac{1-\gamma}{\gamma}} < \cdots < a_N(t)^{\frac{1-\gamma}{\gamma}} < b_N(t)^{\frac{1-\gamma}{\gamma}}. $$
Thus $(V_{A}^{\nu})''(t) \leq 0$. Hence the function $V_{A}^{\nu}$ is concave on $[s_j, s_{j+1})$. Using the same argument as in Case 1, we conclude that the function $V_{A}^{\nu}$ is $s$-concave on $[t_1 + \eps, t_2]$.

If $m \gamma < 0$, then one has 
$$ a_1(t)^{\frac{1-\gamma}{\gamma}} > b_1(t)^{\frac{1-\gamma}{\gamma}} > \cdots > a_N(t)^{\frac{1-\gamma}{\gamma}} > b_N(t)^{\frac{1-\gamma}{\gamma}} \geq 0. $$
We deduce that
\begin{eqnarray*}
\left( \sum_{i=1}^{N} \left( a_i(t)^{\frac{1+\gamma}{\gamma}} - b_i(t)^{\frac{1+\gamma}{\gamma}} \right) \right) \left( \sum_{i=1}^{N-1} \left( a_i(t)^{\frac{1-\gamma}{\gamma}} - b_i(t)^{\frac{1-\gamma}{\gamma}} \right) + a_N(t)^{\frac{1-\gamma}{\gamma}} \right) \\ \leq a_1(t)^{\frac{1+\gamma}{\gamma}} a_1(t)^{\frac{1-\gamma}{\gamma}} \leq \left( \sum_{i=1}^{N-1} \left( b_i(t)^{\frac{1}{\gamma}} + a_i(t)^{\frac{1}{\gamma}} \right) + a_N(t)^{\frac{1}{\gamma}} \right)^2.
\end{eqnarray*}
Hence the function $V_{A}^{\nu}$ is $s$-concave on $[s_j, s_{j+1})$. We conclude that $V_A^{\nu}$ is $s$-concave on $[t_1 + \eps, t_2]$. \\

Hence we get that $\Phi_{\eps}(\nu) \leq 0$ and it follows that $V_A^{\mu}$ is $s$-concave on $(0, t_0)$. We have seen that $V_A^{\mu}$ is $s$-concave on $[t_0, + \infty)$ and using Lemma~\ref{raccord} we conclude that $V_A^{\mu}$ is $s$-concave on $(0, + \infty)$. Finally, from the non-decreasing property of $V_A^{\mu}$, it follows that $V_A^{\mu}$ is $s$-concave on $\R_+$.
\end{proof}

\noindent
{\bf Remark.}
The result clearly holds if we replace the interval $[-1,1]$ by any symmetric interval. However it is not necessarily true for an arbitrary interval. For example, let $0 < s \leq \frac{1}{2}$ and take $B = [0, 1]$, $A = [0,1] \cup [2,3]$ and $\, \de \mu(x) = x^{\frac{1}{\gamma}} 1_{[0,3]}(x) \, \de x$, with $\gamma = \frac{s}{1-s}$.
Then $\mu$ is an $s$-concave measure. For $t \in [0,\frac{1}{2})$ we get
$$ V_A^{\mu}(t)=\mu(A + tB) = \frac{\gamma}{\gamma + 1} \left( (1+t)^{\frac{\gamma + 1}{\gamma}} + 3^{\frac{\gamma + 1}{\gamma}} - 2^{\frac{\gamma + 1}{\gamma}} \right). $$
Thus,
$$ V_A^{\mu}(0)(V_A^{\mu})''(0) - (1-s)(V_A^{\mu})'(0)^2 = \frac{1}{\gamma + 1} \left( 3^{\frac{\gamma + 1}{\gamma}} - 2^{\frac{\gamma + 1}{\gamma}} \right) > 0. $$
Hence $V_A^{\mu}$ is not $s$-concave on $\R_+$. For $s=0$, the same example works. For $s < 0$, one can take $B = [-1, 0]$, $A = [0,1] \cup [2,3]$ and $\, \de \mu(x) = x^{\frac{1}{\gamma}} 1_{[\alpha,3]}(x) \, \de x$, with $\gamma = \frac{s}{1-s}$ and $\alpha$ sufficiently small. \\

The localization theorem (Theorem~\ref{FG}) holds only for $s \leq \frac{1}{2}$. Thus this theorem could not be used to examine the case of $s$-concave measures, with $s > \frac{1}{2}$. For $s \in (\frac{1}{2},1)$, the answer to Problem~1 is negative in general as shown in Proposition~\ref{gap} below, but under specific conditions regarding the support of the measure, we can show that a stronger positive answer holds. First, let us show that for $s \in (\frac{1}{2},1)$ the answer to Problem~1 is negative in dimension 1.

\begin{prop}\label{gap}

Let $s \in (\frac{1}{2}, 1)$ and let $\gamma = \frac{s}{1-s}$. Let $\beta= 10 (1-2^{\frac{1-\gamma}{\gamma}})^{-1}$ and let $\mu$ be a measure such that $\de \mu(x) = x^{\frac{1}{\gamma}}1_{[0,\beta]}(x) \de x$. Let us set $A=[0,1] \cup [2,\beta]$. Then the function $t \mapsto V_A^{\mu}(t) = \mu(A + t[-1,1])$ is not $s$-concave on $\R_+$.

\end{prop}

\begin{proof}
For every $t \in [0,\frac{1}{2})$,
\begin{eqnarray*}
V_A^{\mu}(t) & = & \frac{\gamma}{\gamma + 1} \left( (1+t)^{\frac{1+\gamma}{\gamma}} + \beta^{\frac{1+\gamma}{\gamma}} - (2-t)^{\frac{1+\gamma}{\gamma}} \right), \\ (V_A^{\mu})'(t) & = & (1+t)^{\frac{1}{\gamma}} + (2-t)^{\frac{1}{\gamma}}, \\ (V_A^{\mu})''(t) & = & \frac{1}{\gamma} \left( (1+t)^{\frac{1-\gamma}{\gamma}} - (2-t)^{\frac{1-\gamma}{\gamma}} \right).
\end{eqnarray*}
Hence,
$$ V_A^{\mu}(0) (V_A^{\mu})''(0) - (1-s) (V_A^{\mu})'(0)^2 \! = \! \frac{1}{\gamma + 1} \! \left(\! \beta^{\frac{1+\gamma}{\gamma}}(1-2^{\frac{1-\gamma}{\gamma}}) - 2^{\frac{1-\gamma}{\gamma}} - 2^{\frac{1+2\gamma}{\gamma}} \right) \!. $$
Since $1-2^{\frac{1-\gamma}{\gamma}} > 0$ and $\beta>\left( (2^{\frac{1+2\gamma}{\gamma}} + 2^{\frac{1-\gamma}{\gamma}}) / (1-2^{\frac{1-\gamma}{\gamma}}) \right)^{\frac{\gamma}{\gamma+1}}$, it follows that $V_A^{\mu}(0) (V_A^{\mu})''(0) - (1-s) (V_A^{\mu})'(0)^2 > 0$. We conclude that $V_A^{\mu}$ is not $s$-concave on $\R_+$.
\end{proof}

We denote by $\supp(\mu)$ the support of $\mu$ and by $\dist(A, \supp(\mu)^c)$ the distance between $A$ and the complement of the support of $\mu$. When the support of $\mu$ is $\R$, the distance will be equal to $+ \infty$.

\begin{prop}\label{theo:3.7}

Let $s \geq \frac{1}{2}$. Let $\mu$ be an $s$-concave measure on $\R$ which is absolutely continuous with respect to the Lebesgue measure. Let $A$ be a compact subset of $\R$ such that $\dist(A, \supp(\mu)^c))>0$. Then the function $t \mapsto V_A^{\mu}(t) = \mu(A + t[-1,1])$ is concave on $[0, \dist(A, \supp(\mu)^c))$.

\end{prop}

\begin{proof}
First, we assume that $s = \frac{1}{2}$. Hence $\mu$ admits a $1$-concave density function denoted by $\psi$. Notice that $V_A^{\mu}$ is differentiable outside a finite number of points $t_0, \dots, t_N$ and that for every $t \in [0, \dist(A, \supp(\mu)^c)) \setminus \{t_0, \cdots, t_N\}$ one has
$$ V_A^{\mu}(t) = \sum_{i = 1}^N \int_{a_i - t}^{b_i + t} \psi(x) \, \de x, \quad (V_A^{\mu})'(t) = \sum_{i=1}^N \left( \psi(b_i+t) + \psi(a_i-t) \right). $$
Since $\psi$ is concave, it follows that for every $i \in \{1,\cdots,N\}$, the function $t \mapsto \psi(b_i+t) + \psi(a_i-t)$ is non-increasing. Thus $(V_A^{\mu})'$ is piecewise non-increasing. We conclude that $V_A^{\mu}$ is piecewise concave on $[0, \dist(A, \supp(\mu)^c))$. From Lemma~\ref{raccord}, we deduce that $V_A^{\mu}$ is concave on $[0, \dist(A, \supp(\mu)^c))$.

Finally, if $\mu$ is $s$-concave with $s \geq \frac{1}{2}$, then $\mu$ is $\frac{1}{2}$-concave and we conclude from the first part of the proof that $V_A^{\mu}$ is concave on $[0, \dist(A, \supp(\mu)^c))$.
\end{proof}

Now we turn to the study of Problem~1 in dimension $n \geq 2$. It was shown in \cite{nous} that the Costa-Cover conjecture is false in dimension $n \geq 2$, and thus the answer to Problem~1 is negative in general. Let us recall the counterexample.

Let $n \geq 2$. Let us set $A=B_2^n \cup \{ 2 e_1\}$ and let us denote $V_A(t) = |A + t B_2^n|$. Then for every $t \in [0,\frac{1}{2})$, one has
$$ V_A(t) = |B_2^n \cup \{ 2 e_1\} + tB_2^n| = |B_2^n + tB_2^n| + |t B_2^n| = |B_2^n|((1+t)^n+t^n). $$
Since the $\frac{1}{n}$-power of this function is not concave (it is strictly convex), $V_A$ is not $\frac{1}{n}$-concave on $\R_+$ for $n \geq 2$. \\

It could appear surprising that Problem~1 has a negative answer in dimension $n \geq 2$ since we proved that this problem has a positive answer in dimension~$1$ with the localization theorem. This localization technique is usually used to reduce inequalities for general convex measures in dimension~$n$ to inequalities for measures for which the support is a segment, thus the problem becomes 1-dimensional (see $e.g.$ \cite{Fradelizi1} and references therein). Let us explain why one cannot reduce Problem~1 to the dimension~$1$. The reduction done in dimension~$1$ with localization works the same way in dimension~$n$ and we get the following equivalence for every compact set $A$ of $\R^n$:\\
i) $V_A^{\mu}$ is $s$-concave for every measure $\mu$ $s$-concave.\\
ii) $V_A^{\nu}$ is $s$-concave for every measure $\nu$ $s$-affine on a segment $[\alpha,\beta]$. \\
However, ii) is not true in dimension $n \geq 2$. Indeed, we can construct an explicit counterexample to show that the function $t \mapsto |(A + tB_2^n) \cap [\alpha,\beta]|_1$ is not continuous everywhere inside its support and hence this function is not $s$-concave. For example, consider $A = \{(0,0)\} \cup \{(3,0)\} \cup \{(x,1) ; x \in [1,2]\}$ and $[\alpha,\beta] = \{(x,0); x \in [0,3]\}$. \\

In dimension 2, it was shown in \cite{nous} that if $A$ is a connected subset of $\R^2$ then the function $t \mapsto |A + t B_2^2|$ is $\frac{1}{2}$-concave on $\R_+$. However, the following proposition shows that this property fails in general if we replace the Lebesgue measure by an arbitrary $s$-concave measure.

\begin{prop}\label{exemple2}

In dimension 2, there exists a connected set $A$ and a $\frac{1}{2}$-concave measure $\mu$ such that $t \mapsto \mu(A+tB_2^2)$ is not $\frac{1}{2}$-concave on $\R_+$.

\end{prop}

\begin{proof}
We set $\de \mu(x) = 1_{B_1^2}(x) \de x$, where $B_1^2$ denotes the unit ball for the $\ell_1^2$ norm. Hence $\mu$ is $\frac{1}{2}$-concave. We construct the points $B=(-1,0)$, $C=(-0.5,-0.5)$, $D=(0.5,0.5)$, $E=(0,1)$, $F=(-2,0)$, $G=(0,-2)$, $H=(0,-1)$, $I=(2,0)$, $J=(1,0)$. We set
$$ A = \conv(BCDE) \cup [FB] \cup [FG] \cup [GH] \cup [GI] \cup [IJ]. $$
Then $A$ is connected and for every $t \in [0, \frac{1}{8}]$, we get
$$ V_A^{\mu}(t) = \mu(A + t B_2^2) = \frac{\sqrt{2}}{2} + \sqrt{2} t + \frac{\pi}{2} t^2. $$
It follows that $\left( \sqrt{V_A^{\mu}} \right)''(0) > 0$. We conclude that $t \mapsto \mu(A + t B_2^2)$ is not $\frac{1}{2}$-concave on $\mathbb{R}_+$.
\end{proof}

\noindent
{\bf Remark.}
Notice that we can adapt the counterexample of Proposition~\ref{exemple2} to show that there exists an $s$-concave measure $\mu$ on $\R^n$, $n \geq 2$, such that for every $r \in (-\infty, s)$ there exists a compact connected set $A \subset \R^n$ such that $t \mapsto \mu(A+tB_2^n)$ is not $r$-concave on $\R_+$. \\

\noindent
{\bf Question.}
Does there exist $s \leq \frac{1}{2}$ such that for every $s$-concave measure $\mu$ on $\R^2$ and for every compact set $A \subset \R^2$ the function $t \mapsto \mu(A+tB_2^2)$ is $s$-concave for the $t$'s so that the set $\supp(\mu) \cap (A+tB_2^2)$ is connected? 

In particular, for every convex subset $K$ of $\R^2$ and for every compact set $A \subset \R^2$, is the function $t \mapsto |(A+tB_2^2) \cap K|$ $\frac{1}{2}$-concave for the $t$'s so that the set $K \cap (A+tB_2^2)$ is connected?

\section{Functional version}

In this section, we give a partial answer to Problem~2.

As in section 3 where Problem~1 was solved for convex sets, the next proposition shows that Problem~2 is solved for $\gamma$-concave functions, as is expected.

\begin{prop}\label{theo:5.2}

Let $\gamma \geq -\frac{1}{n}$. Let $f,g : \R^n \to \R_+$ be two $\gamma$-concave functions. Then the function $ t \mapsto \int_{\R^n} h_t^{(\gamma)}(z) \, \de z $ is $\frac{\gamma}{1+\gamma n}$-concave on $\R_+$, where $$ h_t^{(\gamma)}(z) = \underset{f(x)>0;~g(y)>0}{\sup_{z=x+ty}}(f(x)^{\gamma} + tg(y)^{\gamma})^{\frac{1}{\gamma}} \quad \mathrm{and} \quad h_t^{(0)}(z) = \sup_{z=x+ty}f(x)g(y)^t. $$

\end{prop}

\begin{proof}
We examine the case $\gamma \neq 0$, the case $\gamma = 0$ can be proved with the same argument. 

For convenience, let us denote $h_t = h_t^{(\gamma)}$. Let $\lambda \in [0,1]$ and let $t_1,t_2 \in \R_+$. We want to show that
$$ \int_{\R^n} h_{(1-\lambda)t_1 + \lambda t_2} \geq \left( (1-\lambda)\left( \int_{\R^n} h_{t_1} \right)^{\frac{\gamma}{1+\gamma n}} + \lambda \left( \int_{\R^n} h_{t_2} \right)^{\frac{\gamma}{1+\gamma n}} \right)^{\frac{1+\gamma n}{\gamma}}. $$
From the Borell-Brascamp-Lieb inequality \cite{Borell}, \cite{Brascamp} (dimensional Pr\'ekopa's inequality), it is sufficient to show that for every $z_1,z_2 \in \R^n$ one has
$$ h_{(1-\lambda)t_1 + \lambda t_2}((1-\lambda)z_1 + \lambda z_2) \geq \left( (1-\lambda)h_{t_1}(z_1)^{\gamma} + \lambda h_{t_2}(z_2)^{\gamma} \right)^{\frac{1}{\gamma}}. $$
Let $z_1,z_2 \in \R^n$. Let $x_1, x_2 \in \R^n$ such that
$$ \forall i \in \{1,2\}, \, h_{t_i}(z_i) = \left( f(x_i)^{\gamma} + t_i g \left( \frac{z_i-x_i}{t_i} \right)^{\gamma} \right)^{\frac{1}{\gamma}}. $$
Let us denote $h = h_{(1-\lambda)t_1 + \lambda t_2}((1-\lambda)z_1 + \lambda z_2)$ and $t=(1-\lambda)t_1 + \lambda t_2$. We get,
\begin{eqnarray*}
h \!\!\! & = & \!\!\! \sup_{x \in \R^n} \left(f(x)^{\gamma} + t g \!\left( \frac{(1-\lambda)z_1 + \lambda z_2 - x}{t} \right)^{\gamma} \right)^{\frac{1}{\gamma}} \\ & \geq & \!\!\! \left( f((1-\lambda)x_1 + \lambda x_2)^{\gamma} + t g \!\left(\! \frac{(1-\lambda)z_1 + \lambda z_2 - ((1-\lambda)x_1 + \lambda x_2)}{t} \!\!\right)^{\gamma} \right)^{\frac{1}{\gamma}} \\ & \geq & \!\!\! \left(\! (1-\lambda)f(x_1)^{\gamma} + \! \lambda f(x_2)^{\gamma} + (1-\lambda) t_1 g \!\left( \frac{z_1-x_1}{t_1} \right)^{\!\! \gamma} \!\!\! + \! \lambda t_2 g \!\left( \frac{z_2-x_2}{t_2} \right)^{\!\!\gamma} \right)^{\!\!\frac{1}{\gamma}} \\ & = & \!\!\! \left( (1-\lambda) h_{t_1}(z_1)^{\gamma} + \lambda h_{t_2}(z_2)^{\gamma} \right)^{\frac{1}{\gamma}}.
\end{eqnarray*}
\end{proof}

As a consequence of the H\"older inequality, if $f: \R^n \to \R_+$ is $\beta$-concave and if $g: \R^n \to \R_+$ is $\gamma$-concave, then $fg$ is $\alpha$-concave for every $\alpha, \beta, \gamma \in \R \cup \{ +\infty \}$ such that $\beta + \gamma \geq 0$ and $\frac{1}{\beta} + \frac{1}{\gamma} = \frac{1}{\alpha}$ (see $e.g.$ \cite{HLP}). A generalized form of Proposition~\ref{theo:5.2} follows:

\begin{prop}

Let $\gamma \geq -\frac{1}{n}$. If a measure $\mu$ has a $\beta$-concave density, with $\beta \geq -\gamma$, and if $f,g : \R^n \to \R_+$ are two $\gamma$-concave functions, then $t \mapsto \int_{\R^n} h_t^{(\gamma)}(z) \, \de \mu(z)$ is $\frac{\alpha}{1+\alpha n}$-concave on $\R_+$, where $\frac{1}{\beta} + \frac{1}{\gamma} = \frac{1}{\alpha}$.

\end{prop}

In Problem~2, for $\gamma < 0$ the function $V=g^{\gamma}$ is convex by assumption and one can naturally connect the function $h_t^{(\gamma)}$ with the Hopf-Lax solution of the Hamilton-Jacobi equation:
$$ h_t^{(\gamma)}(z) = \sup_{x \in \R^n} \left(f(x)^{\gamma} + t V \left(\frac{z-x}{t} \right) \right)^{\frac{1}{\gamma}} = \left( Q_t^{(V)} f^{\gamma}(z) \right)^{\frac{1}{\gamma}}, $$
where for arbitrary convex function $V$ and for arbitrary function $u$,
$$ Q_t^{(V)} u(z) = \inf_{x\in \R^n} \left( u(x) + t V \left(\frac{z-x}{t} \right) \right). $$

Let us recall basic properties of the Hopf-Lax solution of the Hamilton-Jacobi equation. For a convex function $V$ such that $\lim_{|z| \to + \infty} V(z) / |z| = + \infty$ and for a Lipschitz function $u$, it is known that $Q_t^{(V)} u$ is the solution, called {\it Hopf-Lax solution}, of the following partial differential equation, called {\it Hamilton-Jacobi equation} (see $e.g.$ \cite{Evans}).
\[ 
\left\{
  \begin{array}{ll}
     \frac{\partial}{\partial t} h(t,z) + V^*(\nabla_z h(t,z)) = 0 & \mbox{ on $(0, +\infty) \times \R^n$}\\
     h(t,z)=u(z) & \mbox{ on $\{t=0\} \times \R^n$}
  \end{array},
\right. \]
where $V^*$ is the Legendre transform of $V$ defined on $\R^n$ by 
$$ V^*(y) = \sup_{x \in \R^n} (\langle x,y \rangle - V(x)). $$
It is shown in \cite{Evans} that if $u$ is Lipschitz on $\R^n$ then $Q_t^{(V)} u$ is Lipschitz on $[0, + \infty) \times \R^n$. However, for an arbitrary convex function $V$, $t \mapsto Q_t^{(V)} u$ is not necessarily continuous at $0$.

\begin{proof}[Proof of Theorem 2.]
We examine the case $\gamma \neq 0$, the case $\gamma = 0$ can be proved with the same argument. Since $\gamma \in [-1, 0)$, one has
$$ h_t^{(\gamma)}(z) = \sup_{x \in \R} \left(f(x)^{\gamma} + t V \left(\frac{z-x}{t} \right) \right)^{\frac{1}{\gamma}} = \left( Q_t^{(V)} f^{\gamma}(z) \right)^{\frac{1}{\gamma}}. $$
We denote for $t \in \R_+$, 
$$ F(t) = \int_{\R} h_t^{(\gamma)}(z) \, \de z = \int_{\R} \left( Q_t^{(V)}f^{\gamma}(z) \right)^{\frac{1}{\gamma}} \, \de z. $$
For $p=1$, the function $F$ is constant. We then consider $p>1$. For convenience, we set $\phi(t,z) = Q_t^{(V)} f^{\gamma}(z)$ and $\phi'=\frac{\partial \phi}{\partial z}$. By assumption, there exists $t_0>0$ such that for almost every $z \in \R$, $t \mapsto h_t^{(\gamma)}(z)$ is twice differentiable on $[0,t_0]$ with derivatives that are bounded by an integrable function, hence we get for every $t \in [0,t_0]$,
$$ F'(t) = -\frac{1}{\gamma} \int_{\R} V^*\left( \phi' \right) \phi^{\frac{1-\gamma}{\gamma}}, $$
$$ F''(t) = \frac{1}{\gamma} \int_{\R} \phi''  \left( (V^*)'(\phi') \right)^2 \phi^{\frac{1-\gamma}{\gamma}} + \frac{1-\gamma}{\gamma^2} \int_{\R}(V^*(\phi'))^2 \phi^{\frac{1-2\gamma}{\gamma}}. $$
We assumed that $V(u) = \frac{|u|^p}{p}$. Hence $V^*(u) = \frac{|u|^q}{q}$, with $\frac{1}{p} + \frac{1}{q} = 1$. It follows that
$$ F''(t) = \frac{1}{\gamma} \int_{\R} \phi''  (\phi')^{2q-2} \phi^{\frac{1-\gamma}{\gamma}} + \frac{1-\gamma}{\gamma^2} \int_{\R}\frac{(\phi')^{2q}}{q^2} \phi^{\frac{1-2\gamma}{\gamma}}. $$
By assumption, $\lim_{z \to \pm \infty}\frac{\partial h_t^{(\gamma)}}{\partial z} (z) = 0$, thus integration by parts gives
$$ \frac{2q-1}{\gamma} \int_{\R} \phi''  (\phi')^{2q-2} \phi^{\frac{1-\gamma}{\gamma}} = - \frac{1-\gamma}{\gamma^2} \int_{\R} (\phi')^{2q}\phi^{\frac{1-2\gamma}{\gamma}}. $$
Finally,
$$ F''(t) = -\frac{1-\gamma}{\gamma^2}\frac{(q-1)^2}{q^2(2q-1)} \int_{\R} (\phi')^{2q}\phi^{\frac{1-2\gamma}{\gamma}} \leq 0. $$
We conclude that $ t \mapsto \int_{\R} h_t^{(\gamma)}(z) \, \de z $
is concave on $[0,t_0]$.
\end{proof}

\noindent
{\bf Open problem.}
Problem~2 is open in dimension $1$ for arbitrary $\gamma$-concave function $g$.

\section{Links with weighted Brascamp-Lieb-type inequalities}

In this section, we express the $s$-concavity of the function $t \mapsto \int_{\R^n} h_t^{(\gamma)}(z) \, \de z$ in term of a weighted Brascamp-Lieb-type inequality.

\begin{prop}\label{BL}

Let $\gamma \in [-\frac{1}{n}, 0)$ and let $s \in \R$. Let $f: \R^n \to \R_+$ be such that $f^{\gamma}$ is a bounded Lipschitz function. Let $V : \R^n \to \R_+$ be a convex function such that $\lim_{|z| \to +\infty} V(z)/|z| = + \infty$. Let us define for every $z \in \R^n$,
$$ h_t^{(\gamma)}(z) = \underset{f(x)>0;~V(y)>0}{\sup_{z=x+ty}}(f(x)^{\gamma} + tV(y))^{\frac{1}{\gamma}}. $$
If there exists $t_0 > 0$ such that for almost every $z \in \R$, $t \mapsto h_t^{(\gamma)}(z)$ is twice differentiable on $[0,t_0]$ with derivatives that are bounded by an integrable function and if $G \in L^2(\mu)$, where $\de \mu(z) = \left( (Q_t^{(V)}f^{\gamma}(z))^{\frac{1}{\gamma}} / \int (Q_t^{(V)}f^{\gamma})^{\frac{1}{\gamma}} \right) \, \de z$ and $ G = V^*(\nabla_z Q_t^{(V)}f^{\gamma}) / Q_t^{(V)}f^{\gamma}$, then the function $t \mapsto \int_{\R^n} h_t^{(\gamma)}(z) \, \de z$ is $s$-concave on $[0,t_0]$ if and only if
\begin{eqnarray*}
\Var_{\mu} \left( G \right) & \leq & - \frac{\gamma}{1-\gamma} \int \frac{ \langle \, (\Hess\, Q_t^{(V)}f^{\gamma}) (\nabla_z V^*)(\nabla_z Q_t^{(V)}f^{\gamma}),(\nabla_z V^*)(\nabla_z Q_t^{(V)}f^{\gamma}) \, \rangle }{Q_t^{(V)}f^{\gamma}} \, \de \mu \\ & ~ & \qquad \qquad \qquad + \frac{\gamma - s}{1 - \gamma} \left( \int G \, \de \mu \right)^2,
\end{eqnarray*}

\end{prop}

\begin{proof}
Recall that $h_t^{(\gamma)}(z) = \left( Q_t^{(V)} f^{\gamma} (z) \right)^{\frac{1}{\gamma}}$, where for arbitrary convex function $V$ and arbitrary function $u$,
$$ Q_t^{(V)} u(z) = \inf_{x\in \R^n} \left( u(x) + t V \left(\frac{z-x}{t} \right) \right). $$
For convenience, let us denote $\phi = f^{\gamma}$ and $Q_t=Q_t^{(V)}$. We get
$$ \frac{\partial h_t^{(\gamma)}}{\partial t} = -\frac{1}{\gamma} V^*(\nabla_z Q_t \phi) (Q_t \phi)^{\frac{1-\gamma}{\gamma}}, $$
\begin{eqnarray*}
\frac{\partial^2 h_t^{(\gamma)}}{\partial t^2} & = & \frac{1}{\gamma} \langle \, (\Hess\, Q_t \phi) (\nabla_z V^*)(\nabla_z Q_t \phi),(\nabla_z V^*)(\nabla_z Q_t \phi) \, \rangle (Q_t \phi)^{\frac{1-\gamma}{\gamma}} \\ & ~ & \qquad \qquad + \frac{1-\gamma}{\gamma^2} (V^*(\nabla_z Q_t \phi))^2 (Q_t \phi)^{\frac{1-2\gamma}{\gamma}}.
\end{eqnarray*}
Thus the function $F$ is $s$-concave if and only if $ F(t) F''(t) \leq (1-s) F'(t)^2 $ if and only if
\begin{eqnarray*}
\Var_{\mu} \left( G \right) & \leq & - \frac{\gamma}{1-\gamma} \int \frac{ \langle \, (\Hess\, Q_t \phi) (\nabla_z V^*)(\nabla_z Q_t \phi),(\nabla_z V^*)(\nabla_z Q_t \phi) \, \rangle }{Q_t \phi} \, \de \mu \\ & ~ & \qquad \qquad \qquad + \frac{\gamma - s}{1 - \gamma} \left( \int G \, \de \mu \right)^2,
\end{eqnarray*}
where $\de \mu(z) = \left( (Q_t \phi(z))^{\frac{1}{\gamma}} / \int (Q_t \phi)^{\frac{1}{\gamma}} \right) \, \de z$ and $G = \left(V^*(\nabla_z Q_t \phi) / Q_t \phi \right)$.
\end{proof}

\noindent
{\bf Remark.}
For $\gamma = 0$ and $V(u)=\frac{|u|^2}{2}$, one may use the same argument to get that $t \mapsto \int_{\R^n} h_t^{(0)}(z) \, \de z$ is log-concave if and only if
$$ \Var_{\mu}(|\nabla_z Q_t \phi|^2) \leq 4 \int \langle \, (\Hess \, Q_t \phi) \nabla_z Q_t \phi,\nabla_z Q_t \phi \, \rangle \, \de \mu, $$
where $\phi = - \log f$ and $\de \mu(z) = \left( \e^{-Q_t \phi(z)} / \int \e^{-Q_t \phi} \right) \, \de z $. \\

\noindent
{\bf Question.}
For which function $u$ does the following inequality holds?
$$ \Var_{\mu}(|\nabla_z u|^2) \leq 4 \int \langle \, (\Hess \, u) \nabla_z u,\nabla_z u \, \rangle \, \de \mu, $$
where $\de \mu(z) = \left( \e^{-u(z)} / \int \e^{-u} \right) \, \de z $. \\

From Proposition~\ref{theo:5.2}, if $f$ is $\gamma$-concave with $\gamma \in [-\frac{1}{n},0)$, then one may apply Proposition~\ref{BL} to get the following weighted Brascamp-Lieb-type inequality by letting $t$ go to $0$:

\begin{coro}

Let $\gamma \in [-\frac{1}{n}, 0)$ and let $s = \frac{\gamma}{1 + \gamma n}$. For every $V,\phi : \R^n \to \R_+$ convex such that $\lim_{|z| \to +\infty} V(z)/z = + \infty$, one has
\begin{eqnarray}\label{poincare}
\Var_{\mu}(G) \leq - \frac{\gamma}{1-\gamma} \int \frac{ \langle \, (\Hess \phi)^{-1} \nabla_z G \phi, \nabla_z G \phi \, \rangle }{\phi} \de \mu + \frac{\gamma - s}{1 - \gamma} \left( \int G \, \de \mu \right)^2,
\end{eqnarray}
where $\de \mu(z) =  \left( \phi^{\frac{1}{\gamma}}(z) / \int \phi^{\frac{1}{\gamma}} \right) \de z$ and $G= V(\nabla_z \phi) / \phi$.

\end{coro}

We have reproved Bobkov-Ledoux's result \cite{Bobkov2} (for a smaller class of functions $G$) who used the same idea since inequality~(\ref{poincare}) is derived from the Borell-Brascamp-Lieb inequality (dimensional Pr\'ekopa's inequality). Bobkov and Ledoux have already seen in \cite{Bobkov1} that one can deduce the classical Brascamp-Lieb inequality from the classical Pr\'ekopa inequality (corresponding to the log-concave case). This idea has been explored by Cordero-Erausquin and Klartag in \cite{Cordero-Klartag} where they showed that the converse is true, $i.e.$ one can derive the Pr\'ekopa inequality from the Brascamp-Lieb inequality. Thereafter, Nguyen \cite{Nguyen} generalized the work by Cordero-Erausquin and Klartag to the case of $s$-concave measures (even for $s \geq 0$) and improved Bobkov-Ledoux's Brascamp-Lieb-Type inequality (inequality~(\ref{poincare})). Recently, Kolesnikov and Milman \cite{KM} generalized the weighted Brascamp-Lieb-type inequalities obtained by Bobkov, Ledoux and Nguyen to the setting of Riemannian manifolds. \\

\noindent
{\bf Acknowledgements.}\\
I would like to thank my advisor Matthieu Fradelizi, Pierre Youssef and the referee for their valuable comments and suggestions.


\begin{thebibliography}{99}
\bibliographystyle{unsrt}
	
	\bibitem{AKM} S. Artstein, B. Klartag, V. Milman, {\it The Santal\'o point of a function and a functional form of
Santal\'o inequality}, Mathematika 51, 33-48 (2004)
	\bibitem{Ball} K. Ball, {\it Isometric problems in $\ell_p$ and sections of convex sets}, PhD Dissertation, Cambridge (1986)
	\bibitem{Bobkov4} S. G. Bobkov, {\it Isoperimetric and analytic inequalities for log-concave probability measures}, Ann. Probab. 27 (1999), no. 4, 1903-1921
	\bibitem{Bobkov3} S. G. Bobkov, I. Gentil, M. Ledoux, {\it Hypercontractivity of Hamilton-Jacobi equations}, J. Math. Pures Appl. (9) 80 (2001), no. 7, 669-696	
	\bibitem{Bobkov1} S. G. Bobkov, M. Ledoux, {\it From Brunn-Minkowski to Brascamp-Lieb and to logarithmic Sobolev inequalities}, Geom. Funct. Anal. 10 (2000), no. 5, 1028-1052
	\bibitem{Bobkov2} S. G. Bobkov, M. Ledoux, {\it Weighted Poincar\'e-type inequalities for Cauchy and other convex measures}, Ann. Probab. 37 (2009), no. 2, 403-427
	\bibitem{Borell1} C. Borell, {\it Convex measures on locally convex spaces}, Ark.Mat. 12 (1974), 239-252
	\bibitem{Borell} C. Borell, {\it Convex set functions in d-space}, Periodica Mathematica Hungarica Vol. 6, 111-136, 1975
	\bibitem{Borell0} C. Borell, {\it The Brunn-Minkowski inequality in Gauss space}, Invent. Math. 30 (1975), no. 2, 207-216
	\bibitem{Brascamp} H. J. Brascamp, E. H. Lieb, {\it On extensions of the Brunn-Minkowski and Pr\'ekopa-Leindler theorems, including inequalities for log concave functions, and with an application to the diffusion equation}, (1976) J. Funct. Anal. 22 366-389.
	\bibitem{Cianchi} A. Cianchi, N. Fusco, F. Maggi, A. Pratelli, {\it On the isoperimetric deficit in Gauss space}, Amer. J. Math. 133 (2011), no. 1, 131-186 
	\bibitem{Cordero-Klartag} D. Cordero-Erausquin, B. Klartag, {\it Interpolations, convexity and geometric inequalities}, Geometric aspects of functional analysis, 151-168, Lecture Notes in Math., 2050, Springer, Heidelberg, 2012
	\bibitem{C} M. Costa, {\it A new entropy power inequality},  IEEE Trans. Inform. Theory 31 (1985) 751-760
	\bibitem{Costa} M. Costa, T. M. Cover, {\it On the similarity of the entropy power inequality and the Brunn-Minkowski inequality}, IEEE Trans. Inform. Theory 30 (1984), no. 6, 837-839
	\bibitem{Castro} Y. De Castro, {\it Quantitative isoperimetric inequalities on the real line}, Ann. Math. Blaise Pascal 18 (2011), no. 2, 251-271
	\bibitem{Evans} L. C. Evans, {\it Partial differential equations}, Graduate Studies in Mathematics, 19. American Mathematical Society, Providence, RI, 1998. xviii+662 pp
	\bibitem{Fradelizi1} M. Fradelizi, {\it Concentration inequalities for $s$-concave measures of dilations of Borel sets and applications}, Electron. J. Probab. 14 (2009), no. 71, 2068-2090
	\bibitem{Fradelizi} M. Fradelizi, O. Gu\'edon, {\it The extreme points of subsets of s-concave probabilities and a geometric localization theorem}, Discrete Comput. Geom. 31 (2004), no. 2, 327-335
	\bibitem{nous} M. Fradelizi, A. Marsiglietti, {\it On the analogue of the concavity of entropy power in the Brunn-Minkowski theory},  Adv. in Appl. Math. 57 (2014), 1-20
	\bibitem{FradeliziM1} M. Fradelizi, M. Meyer, {\it Some functional forms of Blaschke-Santal\'o inequality}, Math. Z. 256 (2007), no. 2, 379-395
	\bibitem{FradeliziM2} M. Fradelizi, M. Meyer, {\it Functional inequalities related to Mahler's conjecture}, Monatsh. Math. 159 (2010), no. 1-2, 13-25
	\bibitem{Gozlan} N. Gozlan, C. Roberto, P-M. Samson, {\it Hamilton Jacobi equations on metric spaces and transport-entropy inequalities}, Rev. Mat. Iberoam. 30 (2014), no. 1, 133-163
	\bibitem{HLP} G. H. Hardy, J. E. Littlewood, G. P\'olya, {\it Inequalities}, Cambridge University Press, Cambridge, 1959
	\bibitem{KLS} R. Kannan, L. Lov\'asz, M. Simonovits, {\it Isoperimetric problems for convex bodies and a localization lemma}, (English summary) Discrete Comput. Geom. 13 (1995), no. 3-4, 541-559
	\bibitem{KM} A. V. Kolesnikov, E. Milman, {\it Poincar\'e and Brunn-Minkowski inequalities on weighted Riemannian manifolds with boundary}, preprint, arXiv:1310.2526 [math.DG]
	\bibitem{Lehec} J. Lehec, {\it A direct proof of the functional Santal\'o inequality}, C. R. Math. Acad. Sci. Paris 347 (2009), no. 1-2, 55-58
	\bibitem{Leindler} L. Leindler, {\it On a certain converse of H\"older's inequality}, II, Acta Sci. Math., 33 (1972), 217-223
	\bibitem{Nguyen} V. H. Nguyen, {\it Dimensional variance inequalities of Brascamp-Lieb type and a local approach to dimensional Pr\'ekopa's theorem}, J. Funct. Anal. 266 (2014), no. 2, 931-955
	\bibitem{Prekopa} A. Pr\'ekopa, {\it On logarithmic concave measures and functions}, Acta Sci. Math. (Szeged) 34 (1973), 335-343
	\bibitem{Sudakov} V. N. Sudakov, B. S. Cirel'son, {\it Extremal properties of half-spaces for spherically invariant measures} (Russian) Problems in the theory of probability distributions, II. Zap. Nau\v{c}n. Sem. Leningrad. Otdel. Mat. Inst. Steklov. (LOMI) 41 (1974), 14-24, 165
	\bibitem{V} C. Villani, {\it A short proof of the 'concavity of entropy power'}, IEEE Trans. Inform. Theory 46 (2000) no. 4, 1695-1696
	
\end{thebibliography}
\end{document}